\documentclass[12pt]{amsart}
\usepackage{amsmath,amsfonts,amsthm,amsopn,cite,mathrsfs}
\usepackage{epsfig,verbatim}
\usepackage{subfigure}

\newtheorem{tm}{Theorem}[section]
\newtheorem{lemma}[tm]{Lemma}
\newtheorem{prop}[tm]{Proposition}
\newtheorem{cor}[tm]{Corollary}

\newcommand{\rem}{\noindent\textsl{REMARK:}}

\newcommand{\field}[1]{\mathbb{#1}}
\newcommand{\bR}{\field{R}}        
\newcommand{\bN}{\field{N}}        
\newcommand{\bZ}{\field{Z}}        
\newcommand{\bC}{\field{C}}        
\newcommand{\bS}{\field{S}}        %
\newcommand{\bP}{\field{P}}        %
\newcommand{\E}{\field{E}\,}        %

 \def\cF{\mathcal{F}}              

 \def\cB{\mathcal{B}}
 \def\cE{\mathcal{E}}

 \def\cA{\mathcal{A}}

 \def\cC{\mathcal{C}}
 
 \def\cO{\mathcal{O}}

 \def\cX{\mathcal{X}}
 \def\cY{\mathcal{Y}}

\def\eps{\epsilon}

\def\rd{\bR^d}
\def\cd{\bC^d}

\def\lrd{L^2(\rd)}

\def\zd{\bZ^d}

\def\intrd{\int_{\rd}}

\def\inv{^{-1}}

\def\wh{\widehat}

\def\E{{{\Bbb E}\,}}
\def\P{{\Bbb P}}

\newcommand{\norm}[1]{\lVert#1\rVert}

\DeclareMathOperator*{\supp}{supp}

 \newcommand{\tf}{time-frequency}

\newcommand{\fif}{if and only if}

\def\Var{{\mathop {{\rm Var\, }}}}

 \theoremstyle{definition}
 \newtheorem{definition}{Definition}

\def\R{\mathbb{R}}

\def\lam{{\lambda}}

\def\brd{\cB (R,\delta )}

\def\vf{\varphi}
\newcommand{\yy}{Y_j(f)}
\newcommand{\jr}{j=1, \dots , r}

\newcommand{\bdl}{bandlimited}
\begin{document}

\begin{abstract}
We consider the problem of random sampling for bandlimited
functions. When  can a bandlimited function $f$ be recovered from
randomly chosen samples $f(x_j), j\in J\subset  \mathbb{N}$? 
We  estimate the probability that a  sampling inequality of
the form 
\begin{equation*}
A\|f\|_2^2 \leq \sum _{j\in J} |f(x_j)|^2 \leq B \|f\|_2^2  
\end{equation*}
 hold uniformly  for all functions $f\in L^2(\mathbb{R}^d)$
with supp $\hat{f} 
\subseteq [-1/2,1/2]^d$ or for  some subset of \bdl\ functions. 

 In contrast to discrete models, the space of bandlimited functions
 is infinite-dimensional and its  
functions ``live'' on the unbounded set $\mathbb{R}^d$. These facts raise  new
problems and leads to  both negative and positive results. 

(a) With probability one, the  sampling inequality  fails  for any
reasonable definition of a random set on 
$\mathbb{R}^d$, e.g., for spatial Poisson processes or  uniform distribution over
disjoint cubes. 

(b)  With overwhelming probability,  the sampling inequality
 holds  
  for certain compact subsets of the space of
bandlimited functions and for  sufficiently large sampling size. 
\end{abstract}

\title[Random Sampling of Bandlimited Functions]{Random Sampling of
Bandlimited Functions} 
\author{Richard F. Bass}
\address{Department of Mathematics, 
The University of Connecticut,
Storrs, CT  06269-3009, USA}
\author{Karlheinz Gr\"ochenig}
\address{Faculty of Mathematics \\
University of Vienna \\
Nordbergstrasse 15 \\
A-1090 Vienna, Austria}
\email{bass@math.uconn.edu,karlheinz.groechenig@univie.ac.at}
\date{}
\thanks{R.B. was partially supported by NSF Grant DMS0601783}
\thanks{K.G. was supported by the  Marie-Curie Excellence Grant
  MEXT-CT-2004-517154} 
\maketitle

\section{Introduction}

The sampling problem asks for the  reconstruction  or approximation of
a function $f$ from its 
sampled values   $  \{f(x_j): j\in J\}$ on some  set 
$\cX = \{x_j \} \subseteq \rd $. In other words, one wants to recover
 $f$ from given samples $f(x_j)$. 
This is a many-faceted problem and spreads over many areas of
mathematics, engineering, and data processing.

We will impose the standard hypothesis that $f$ is bandlimited. In
signal processing  this is a realistic assumption, because it
amounts to assuming a maximum frequency. The assumption is also
relevant in 
complex analysis because a bandlimited function is just  the
restriction of an entire function of 
exponential growth from  $\bC ^d$ to $\rd $. 
 The space of \bdl\ functions is defined to be 
$$
\cB  = \{ f \in \lrd : \supp \hat{f} \subseteq [-1/2,1/2]^d\} \, ,
$$
where we have normalized the spectrum to be the unit cube and the
Fourier transform is normalized as $\hat{f}(\xi ) = \intrd f(x)
e^{-2\pi i x\cdot \xi } \, dx $. 

 The principal   goal is to  establish a  \emph{sampling inequality} of the form
\begin{equation}
  \label{eq:1}
A\|f\|_2^2 \leq \sum _{j} |f(x_j)|^2 \leq B
\|f\|_2^2 \qquad \forall f \in \cB \, .
\end{equation}
A set $\{x_j : j\in J\}\subseteq \rd $ satisfying the sampling
inequality \eqref{eq:1} is called a set of stable sampling or simply a
set of sampling~\cite{landau67}. 
Once a sampling inequality is established,  every $f
\in \cB $ is uniquely determined by its samples on $\cX $ and depends
continuously on these samples. 

Bandlimited functions in dimension $d=1$ and $d>1$ differ in a
fundamental way because of the nature of their zeros.  In dimension
$d=1$ the zeros of an 
entire  function are always discrete, and there is a precise
connection between the possible  density of zeros and the growth of
$f$~\cite{boas,levin,seip04}. By contrast, in higher dimensions, the
zero sets are analytic manifolds, and  standard complex variable
techniques do no longer apply. As a consequence,     almost everything is
known about the sampling of 
bandlimited functions   in dimension $d=1$, but  only a  few results
are known in higher 
dimensions, most notably a strong result of Beurling~\cite{beurling66}. 

The difficulties of the sampling of multivariate functions have
motivated us to turn to probabilistic techniques and to  study
\emph{random sampling}. In this approach the sampling set $\cX $ is a
sequence of random variables $x_j = x_j(\omega )$ on some probability
space $(\Omega, \cF, \P)$ and  taking values in
$\rd $. The sampling inequality~\eqref{eq:1} defines an event on $\Omega
$, and 
 the goal is to estimate the probability that a random set is a set of
 sampling. 

This point of view has worked successfully in  our previous work~\cite{BG04}
where we  have studied the random sampling of
multivariate trigonometric polynomials. We were able to  show that some popular
numerical algorithms~\cite{PST01} work with ``overwhelming''
probability. 
In a similar spirit,   Cand\`es, Romberg, and Tao~\cite{CT06,CRT06}
have recently
investigated sparse trigonometric polynomials and their reconstruction
from a few random samples. 
 The more general context of mathematical learning theory has been
 studied by Cucker, Poggio, Smale, and
 Zhou~\cite{smale1,smale2,SZ04}. In \cite{SZ04} sampling in general
 reproducing kernel Hilbert spaces was studied under the assumption of
 ``rich data.'' This amounts to assuming the validity of a sampling
 inequality. By contrast,  our
interest is  to establish the probability that this basic assumption
holds.  
 The common technical point in  these
approaches \cite{BG04,CT06,smale2} is  the estimate of  entropy and
covering numbers and a metric entropy argument.


The first  contributions to random sampling of \bdl\ functions were
perturbation results in dimension~$d=1$. Seip and
Ulanovsky~\cite{SU96}  investigated  random 
perturbations of regular sampling  $\{j + \delta _j  : j \in \bZ  \}$, where $\delta _j $ is a
sequence of i.i.d. random variables.  
Chistyakov, Lyubarskii, Pastur~\cite{CL97,CLP01}  studied the more general problem
of  perturbation of arbitrary Riesz bases of exponentials. These
contributions are based on the precise characterization of sampling
sets in dimension $d=1$~\cite{seip04}, and the proofs proceed by estimating the
probability that a deterministic condition is satisfied.

For random sampling of \bdl\ functions of several variables new types
of problems arise.

(a) One cannot fall back on deterministic results in higher
dimensions, because sampling theory is not nearly as  developed as in
dimension $d=1$. In fact, this is the very reason why we aim for
purely probabilistic results. 

(b) The space of \bdl\ functions $\cB$ is infinite-dimensional --- in
contrast to trigonometric polynomials of given degree or
sparsity. Thus random matrix techniques as used in~\cite{GPR07,MP06}
are not applicable. 

(c) The configuration space $\rd $ is non-compact and unbounded ---
again  in contrast to trigonometric polynomials that ``live'' on the
torus $[0,1]^d$. 
 This raises the question of how to model a sequence of random points
 in $\rd $.
On a compact set of positive (Lebesgue) measure the
natural notion is that of an independent identically distributed
(i.i.d.) sequence of points with uniform distribution. On $\rd $ there
are several natural choices. We will consider  two such choices: uniform
distributions on disjoint cubes and  spatial Poisson processes. 

%

We will prove that for these two   concepts of ``randomly
distributed points on $\rd $'' the  sampling inequality \eqref{eq:1} must
fail almost surely (Propositions~\ref{RBP1} and~\ref{RBP2}). 
These results come as a surprise to the analyst, but are perhaps more
natural for the probabilist. The reasons for the failure of a sampling
inequality are either the zeros of entire functions or large holes in
the sampling set. In the model of uniform distribution over disjoint
cubes, many  samples may be near the zeros of a bandlimited function 
with positive probability. In other words, the lower bound in
\eqref{eq:1} is small. 

In the other model (spatial Poisson process) we show that, with
positive probability, there are  large  holes   in the
sampling set, which again implies a small lower bound in
\eqref{eq:1}. 

To obtain insight into the formulation of positive results, we argue
in a practical manner. Realistically one can sample $f$ only on a
bounded set; furthermore, every bandlimited function vanishes at
infinity, thus samples far out do not contribute anything significant
to a sampling inequality. We can learn about $f\in \cB $ only if the
samples are taken in the ``essential support'' of $f$, i.e., the set
where most of the $L^2$-norm is localized. Thus we will study the
subset $$\brd =\left\{ f \in \cB : 
    \,\, \int _{[-R/2,R/2]^d} |f(x)|^2 \, dx \geq     (1-\delta )
    \|f\|_2^2 \right\}$$ of 
    bandlimited functions. This subset is compact in $\cB $ and thus
   somewhat resembles a finite-dimensional subspace. Since $f\in \brd
   $ is small outside the cube $[-R/2,R/2]^d$, it should suffice to
   sample $f$ on the relevant cube. In this way, we are back to a
   compact configuration space and an almost finite-dimensional
   function space. Our main result (Theorem~\ref{positive}) is a
   restricted sampling inequality 
   for the subset $\brd $. The proof is a combination of analytic and
   probabilistic    techniques. On the one hand, we will use detailed
   properties about the spectrum of time-limiting operators on \bdl\
   functions by Widom~\cite{widom64}, on the other hand, the metric
   entropy method (see e.g., ~\cite{Du73}). 

The paper is organized as follows. In Section~2 we discuss two natural
models
for random sequences in $\rd $ and show that, with probability one,
they fail to produce sets of stable sampling. In Section~3 we restrict
the attention to  a subset of \bdl\ functions and show that on this
subset a sampling inequality holds with overwhelming probability. The
proof of this  result is contained in Section~4.  In
order to set up the metric entropy method, we
discuss the spectrum of time-frequency limiting operators and covering
numbers. 
We mention that we use two distinct inequalities of Bernstein, one
from Fourier analysis bounding the $L^\infty$ norm of the derivative,
and the other from probability giving estimates for the sums
of independent random variables.

\textbf{Acknowledgement.} We would like to thank the anonymous referee
for his useful comments and pointing out an embarrassing error in the
first version. 

\section{Negative Results}\label{S2}

In  the case of multivariate trigonometric polynomials, 
we showed that if one chose points independently and  uniformly distributed
over the state space, then one could recover the trigonometric polynomial
exactly 
provided only that one had at least as many sample points as the
dimension~\cite[Thm~3.2]{BG04}. 

We first show that this is far from the case for bandlimited functions.
The difficulty is that the state space is not compact.

We first recall a fundamental necessary condition of Landau for a set
of sampling. 
 Let 
 \begin{equation}
   \label{eq:j1}
   D^- (\cX ) = \lim _{R\to \infty } \min _{y\in \rd }
   \frac{\mathrm{card}\, \cX \cap (y+[0,R]^d)}{R^d}
 \end{equation}
be the (lower) Beurling density of a set $\cX\subseteq \rd $. 

\begin{prop}\label{hole}  
Assume that  $\cX = \{x_j\}$ is a set of stable sampling for $\cB
$. Then $\cX  $  must have the following properties: \\
(i) $D^-(\cX ) \geq 1$, in particular there is $R>0$   such that every cube of side
length $R$ contains a sampling point, i.e. 
$\cX \cap (x+ [0,R]^d)
\neq \emptyset $ for all $x\in \rd $.

(ii)  the number of samples in any cube of length $1$ is bounded,
$\max _{y\in \rd } \mathrm{card}\, \cX \cap (y+[0,1]^d )< \infty $. 

A sufficient condition is the following: In dimension $d=1$, if
$D^-(\cX ) >1$ and $\inf _{j\neq k} |x_j - x_k| >0$, then $\cX $ is a
set of sampling.  
\end{prop}

\begin{proof}
(i) is  the result of Landau~\cite{landau67} and have been
re-derived  in ~\cite{GR96a} for discrete sampling sets; the general
case is an  easy extension. 

(ii) is an easy consequence of the finiteness of upper bound $B$ in \eqref{eq:1}. 

 The sufficient condition in dimension $d=1$ is
usually attributed to Beurling and treated in detail by
Seip~\cite{seip04}. 
\end{proof}

 Loosely speaking,  a set of stable sampling must be dense enough and
 cannot have arbitrarily large ``holes''. 


We now consider random sampling sets. Let our probability space be
$(\Omega, \cF, \bP)$ and denote points 
in $\Omega$ by $\omega$.
 When sampling a function $f$
randomly, we consider its samples $f(x_j)$ on a sequence of random
points $x_j = x_j(\omega )$. 
Clearly, a sequence of random points need not have the sufficient
density stated in Proposition~\ref{hole}. However, if the process
is designed to yield  only random sets with $D^-(\cX ) \geq 1$, one
could hope that a generic random set with the necessary density would
be a set of stable sampling. This intuition is completely  false, as we will show
in the next sections. 

\subsection{Uniform distribution on large disjoint cubes}
There are various ways in which one could choose points randomly in 
$\R^d$. 
As a first model we partition $\rd $ into disjoint cubes $k+[0,1]^d,
k\in \zd $, and, in each  cube, we choose $r$ points independently and
uniformly distributed over $k+[0,1]^d$.  Let $\cX$ be the collection
of sample points; $\cX$ is a random set and 
 thus depends on $\omega$. Clearly $D^-(\cX ) = r$ almost surely, so  one
 may expect that $\cX (\omega )$ is a set of stable sampling with high
 probability. 

Our first result says that one cannot obtain a sampling inequality.


\begin{prop}\label{RBP1}
   Let $r\geq 1$ be the number of random samples in each cube
   $k+[0,1]^d$. With probability one the following holds: \\  
 For each $k>0$ there is a function $f_k \in \cB$ such that
$$\sum_{x_i\in \cX(\omega)} |f_k(x_i)|^2\leq \frac{1}{k} \norm{f_k}_2^2.$$
\noindent The function $f_k$ will necessarily depend on $\omega$.

Consequently, a sampling inequality of the form~\eqref{eq:1} is
violated almost surely. 
\end{prop}

\proof
For notational simplicity we give the proof only in dimension $d=1$;
the case of several variables is treated similarly. 

 Let $$g(x)=\frac{\sin(\pi x/2)}{\pi x/2}\, ,$$
and  let
$\psi$ be a nonnegative $C^\infty$ function with support in  
$[-1/4,1/4]$ such that $\psi=1$ on $[-1/8,1/8]$. Let
$\Psi$ be the inverse Fourier transform of $\psi$ 
and define  $F(x)=g(x)\Psi(x)$.

 Since $\psi$ and thus $\Psi$ are  in the Schwartz class, $F$ is in
 $L^2$,  decays
 rapidly,  and there exists a constant $c_1$ such that $|F(x)|\leq c_1/(1+|x|^2)$.
The Fourier transform of $F$ is $\wh g * \psi$, so the support
of $\wh F$ lies in $[-1/2,1/2]$, i.e., $F\in \cB$. Since $F $ is
bounded,   by Bernstein's inequality, $F'$ is also bounded, say,
by $c_2 := \|F'\|_\infty$.

Choose $N$ a large even  integer so that
\begin{equation}\label{RB1.2}
\sum_{|j|\geq N/2} \frac{c_1^2r}{(1+(|j|-1)^2)^2}  <\frac{\norm{F}_2^2}{4k}. 
\end{equation}
Choose $\delta>0$ small so that 
\begin{equation}\label{RB1.3}
2c_2Nr\delta<\frac{\norm{F}_2^2}{2k}. 
\end{equation}
 Let $A_j$ be the event that in the interval $[j,j+1]$ all $r$
points that were chosen randomly lie within $(j, j+\delta)$ if
$j$ is even and within $(j+1-\delta, j+1)$ if $j$ is odd. The events
$A_j$ are independent and the probability of $A_j$ is $\delta^r$.

Let $B=\bigcap_{j=-N}^N A_j$ be the event that the samples in $[-N,N]$
are in a $\delta $-neighborhood of the even integers. By independence,
the probability 
of $B$ is $(\delta^r)^{2N+1}$. If $\omega\in B$, then
using (\ref{RB1.2}) and our bound on $F$,
$$\sum_{x_i\in \cX(\omega)\setminus [-N,N]} |F(x_i)|^2 
\leq \sum_{|j|\geq N/2} \frac{c_1^2 r}{(1+(|j|-1)^2)^2}
<\frac{\norm{F}_2^2}{4k}.$$
By construction, $F(2j)=0$ for $j\in \bZ $, and so using the bound
on $F'$, we have $|F(x)|\leq c_2\delta$ if $|x-2j| \leq
\delta$. Therefore if $\omega\in B$, then 
$$\sum_{x_i\in \cX(\omega)\cap [-N,N]} |F(x_i)|^2
\leq 2c_2Nr\delta<\frac{\norm{F}_2^2}{2k}.$$

Combining, 
if $\omega\in B$, then
\begin{equation}\label{RB1.4}
\sum_{x_i\in \cX(\omega)} |F(x_i)|^2< \frac{\norm{F}_2^2}{k}.
\end{equation}

Now let $C_m =\bigcap_{j=3mN-N}^{3mN+N} A_j$. Clearly the probability
of $C_m$ is the same as the probability of $B$. So
$\sum_{m=1}^\infty \bP(C_m)=\infty$. By independence and the Borel-Cantelli
lemma,
with probability one, $C_m$ occurs for infinitely many $m$. If
$\omega\in C_m$, let $f_k(x)=F(x-3mN)$. Clearly, $f_k \in \cB $ and
the same bounds $c_1$ and $c_2$ 
hold for $f_k$ as for $F$, provided translation is taken into account.
As in (\ref{RB1.4}),
$$\sum_{x_i\in \cX(\omega)} |f_k(x_i)|^2< \frac{\norm{f_k}_2^2}{k}.$$
Thus we have proved that, with probability $1$,   $\cX$  fails to be a set
of stable sampling.  
\qed

\subsection{Spatial Poisson processes}

Another  scheme of choosing points randomly in $\R^d$ is the  spatial
Poisson process $\cX $. This means that 
for some (intensity) function $\lam:\R^d\to [0,\infty)$, for any Borel subset of
$\R^d$, the number of
points in $\cX\cap A$ is  a Poisson random variable with
parameter $\int_A \lam(x) \, dx$.  If $A_1, \ldots, A_n$ are
disjoint sets, then the number of points in $\cX\cap A_i$ are independent
random variables. 

The most natural case is where $\lam(x)$ is a constant,  $\lambda (x)
= \rho $, say. Then the expected Beurling density of $\cX $ is $
\rho $. Again one might think that $\cX $ is a set of stable
sampling with high probability. However,  as in 
Proposition \ref{RBP1}, one cannot get the sampling inequality. In fact,
a stronger result is true. One can choose points at a higher rate
further from the origin and still have the sampling inequality failing.

\begin{prop}\label{RBP2}
  Suppose $\cX$ is a spatial Poisson process
with $\lam(x)=o(1+\log^+(|x|))$. Then with probability one, the
sampling inequality~\eqref{eq:1} fails for \emph{every } subset $\cY $
of $\cX$. 
\end{prop}

\proof Under the hypothesis on $\lambda $, the Beurling density of
$\cX $ may be infinite. In this case, $\cX $ contains too many samples
and the upper bound in the sampling inequality~\eqref{eq:1} will
fail to hold. 
This problem could be fixed by extracting a subset $\cY $ of $ \cX $
that satisfies the necessary conditions of Proposition~\ref{hole}. So
one may still hope that a subsequence   $\cY$ 
may yield a set of stable sampling.

However, we will show that with probability one, for  
each $k>0$ there exists a cube of side length $k$ that contains
no point of $\cX$. Since the maximal hole of a set of stable sampling
is bounded  by Proposition~\ref{hole}(iii), with probability one,  no
subset of $\cX $ can be a set of stable sampling.

The probability that a Poisson random variable with parameter
$\lam$ is equal to zero is $e^{-\lam}$. If $X_i$, $i=1, \ldots, n$, 
are independent Poisson variables with parameters $\lam_i$, resp.,
then
\begin{eqnarray*}
\bP(\hbox{at least one } &X_i & \hbox{ is zero})
=1-\bP(X_1\ne 0, \ldots, X_n\ne 0) \\
&= &1-\prod_{i=1}^n \bP(X_i\ne 0)
=1-\prod_{i=1}^n (1-e^{-\lam_i})\\
&=& 1-\exp\Big(\sum_{i=1}^n \log (1-e^{-\lam_i})\Big)\\
&\geq& 1-\exp\Big(-\sum_{i=1}^n e^{-\lam_i}\Big)\, .
 \end{eqnarray*}

Let $\eps>0$ be chosen later. Choose $m_0>k$ large
so that $\lam(x)\leq \eps \log(|x|)$ if $|x|> m_0$.
For each $m\geq m_0$ we can
find at least $m$ disjoint cubes of side length $k$ lying
in $B(0, 3mk)\setminus B(0,2mk)$; call them $C_{m1}, \ldots,
C_{mm}$. The number of points in 
$\cX$ lying in any one of the $C_{mj}$ is a Poisson random variable with parameter
less than $c_1 \eps (\log m) k^d$. So by the above, 
the probability that at least one of the $C_{mj}$ is empty  
is greater than
$$1-\exp\Big(-me^{-c_1\eps (\log m) k^d}\Big).$$
If we choose $\eps$ so that $c_1\eps k^d\leq 1/2$, then the
above probability is greater than
$$1- \exp\Big(-m^{1/2}\Big),$$
which will be greater than $1/2$ if $m$ is large enough.

Let $D_m$ be the event that at least one of the cubes
$C_{mj}$, $j=1, \ldots, m$, is empty. For $m$ large,
$\bP(D_m)\geq 1/2$, and the $D_m$ are independent. So by
the Borel-Cantelli lemma the event $D_m$ happens for infinitely many
$m$, with probability 1. In particular, there must
be at least one cube of side length $k$ with no points of
$\cX$ in it. 
\qed

On the other hand, the rate of growth $\log^+(|x|)$ is
critical. 
If the intensity function $\lam $ grows faster than 
some multiple of $\log^+ (|x|)$, then the random sequence $\cX $
cannot have large holes.

\begin{prop}\label{RBP3}
  Suppose $\cX$ is a spatial Poisson process with intensity
$\lam(x)\geq c_0(1+\log^+(|x|))$
for all $x$. Fix $\alpha >0$. If $c_0 \geq (d+1)/\alpha ^d$,  then with
probability one,  there exists $R>0$, such that
every cube $\alpha k + [0,\alpha]^d$ for $\alpha |k| \geq R$  contains
at least one point of $\cX $. 
\end{prop}
\proof 
Let $\bS_\alpha$ be the collection of all cubes of the form $\alpha k
+ [0,\alpha ]^d, k\in \zd $. 
We will show that with probability one, all but finitely many cubes 
in $\bS_\alpha$  contain at least one
point of $\cX$. 

Let $C_k$ be the event that the cube $A=\alpha k +[0,\alpha ]^d$
contains no point of $\cX $. 
If $\alpha |k| \geq N$, then $\lambda (x) \geq c_0 \log N$, and thus
$\lambda (A) \geq c_0 \alpha ^d \log N$. 
Thus for $\alpha |k| \geq N$, the probability that this cube is empty
is 
$$
\P (C_k) = e^{-\lambda (A)} \leq e^{-c_0 \alpha
  ^d \log N} = N^{-c_0 \alpha ^d} \, .
$$
If we choose $c_0 \alpha ^d \geq d+1$, then $\sum _{k\in \zd } \P
(C_k) < \infty $. Then
by the Borel-Cantelli lemma, the probability that infinitely
many of the cubes are empty is 0. Therefore
from some $R$ on (depending on $\omega$), all cubes in $\bS_\alpha$ that are at 
least $R$ from the origin are nonempty.
\qed

\section{A Positive Result: Relevant Sampling}

The key to the arguments in Section \ref{S2} was that random sampling sets have either
arbitrarily large holes or can be concentrated  near the zeros
of a bandlimited function. In the former  case we then constructed a
class of 
functions whose main energy is concentrated on the ``hole''; in the
latter case we constructed a class of  functions with prescribed zeros. These
classes then violate the sampling inequality. 

To obtain positive results we change the focus. Since for no
reasonable random sampling set does the norm equivalence~\eqref{eq:1} hold
with positive probability for all band\-limited functions, we will restrict the class of functions
for which we ask~\eqref{eq:1} to hold. The natural idea is to sample a
given $f$ in the region where a significant part of the energy is
located. In other words, we sample in the region of relevant values. 

This idea motivates the following definition. Let $C_R = [-R/2, R/2]^d
$ be the cube of length $R$ centered at the origin. Its volume is
$\mathrm{vol}\, C_R = R^d$. 

\begin{definition}
  Fix a large  number  $R>0$ and a small  $\delta\in (0,1)$.
   Set
  \begin{equation}
    \label{eq:r1a}
    \widetilde{\cB }(R,\delta ) = \left\{ f \in \cB :  \int _{C_R} |f(x)|^2
    \, dx \geq     (1-\delta ) \|f\|_2^2 \right\} 
  \end{equation}
and 
  \begin{equation}
    \label{eq:r1}
    \cB (R,\delta ) = \left\{ f \in \cB : \|f\|_2^2 = 1 \,\, \mathrm{and}
    \,\, \int _{C_R} |f(x)|^2 \, dx \geq     1-\delta   \right\}
  \end{equation}
\end{definition}
Then $\brd $ is the subset of $\cB $ consisting of those \bdl\
functions  whose energy is largely concentrated on the cube $C_R$. Only a
fraction $\delta $ of the total energy is outside this cube. We note
that $\brd $ may be empty when $\delta $ is chosen too small. (For an
estimate of $\delta $ such that $\brd \neq \emptyset$, see Section
3.1) In the following we assume that $\brd $ is non-empty.

 It now
makes sense to sample such $f$ on the cube $C_R$ and to expect that these
samples are  relevant and   capture the main features of
$f$. 

Indeed, we will prove the following result. 

\begin{tm}\label{positive}
  Assume that  $\{x_j : j\in \bN \}$ is a sequence of i.i.d.~random variables
  that are uniformly distributed over the cube $C_R = [-R/2,R/2]^d $
  and $0<\mu <1-\delta $.  Then  there exist $A,B>0$ such that 
the sampling inequality 
\begin{equation}
  \label{eq:cr1}
\frac{r}{R^d}(1- \delta -\mu) \|f\|_2^2 \leq \sum _{j=1}^r |f(x_j)|^2 \leq \frac{r}{R^d}(1+\mu )
\|f\|_2^2 \qquad \forall f \in \widetilde{\cB }(R,\delta )   
\end{equation}
holds  with probability at least
$$ 1- 2Ae^{-B\frac{r}{R^d} \frac{\mu
    ^2}{41+\mu}} \, . $$
The constant  $B $ can be taken to be $B=\frac{\sqrt{2}}{36 }$.  For
large $R$ and sufficiently large sampling size $r $  the constant  $A$  
can be chosen of order $A =  \exp (CR^d )$ with  $C$ 
depending on the dimension $d$.
\end{tm}


\subsection{Discussion and Open Problems}

1. We emphasize  that the exponential probability inequality
holds uniformly for all $f\in \widetilde{\cB }(R,\delta )  $. By
contrast, for fixed $f$ such 
an inequality could be derived much more simply from standard limit
theorems. 

2. Theorem~\ref{positive} is an asymptotic result. It is effective
only for sufficiently large sampling sizes. 
To achieve \eqref{eq:cr1} with a probability exceeding $1-\epsilon $, we
need $2Ae^{-B\frac{r}{R^d} \frac{\mu
    ^2}{41+\mu}} < \epsilon $ or 
\begin{equation}
  \label{eq:cr2}
  r\geq \frac{R^d (41+\mu )}{B \mu ^2} \Big( \log \frac{2}{\epsilon} +
 C R^d \Big) = \cO (R^{2d}) \, .
\end{equation}
Since $\brd $ sits  in a space of approximate  dimension $D=R^d$, we
need $\cO (D^2)$ samples to recover every $f\in \brd $. In finite
dimensional problems, for instance, when sampling trigonometric
polynomials of fixed degree, one can often use random matrix
techniques to show that the effective number of samples is in fact of
the order $\cO (D \log D)$~\cite{GPR07}. It is open whether this
bound is achievable for bandlimited functions in $\brd $.

3. The sampling inequality~\eqref{eq:cr1} states that every $f\in \brd
$ is uniquely determined by a sufficient, but \emph{finite} number of 
samples in $[-R/2,R/2]^d$. This may seem paradoxical at first glance,
because the set of bandlimited functions such that $f(x_j) = 0$ for
$j=1,\dots ,r $,  is an infinite-dimensional subspace of $\cB $. 
However, as we assume that $f$ is essentially supported on the cube 
$[-R/2,R/2]^d$, this means that $f$ must take large values there.  If
$f(x_j)= 0$ for sufficiently many $x_j\in [-R/2,R/2]^d$, then $f$
would oscillate and thus have a large derivative. But this would
contradict the bandlimitedness, which implies that the derivatives of
$f$ are bounded by $\pi$. While Theorem~\ref{positive} is a
probabilistic result, it seems possible to also prove a deterministic
sampling inequality \eqref{eq:cr1} for $\brd $, at least in dimension $d=1$.

4. We emphasize  that $\brd $ is not a subspace. This means
that the frame algorithm (a  linear reconstruction method)
~\cite{DS52} cannot be
used to recover $f$ from its samples. Likewise, the
projection-onto-convex-sets (POCS) method cannot be applied, because $\brd $
is not convex. Although \eqref{eq:cr1} determines each $f\in \brd $
uniquely, currently we do not have an explicit reconstruction
algorithm to recover $f$ from its relevant samples.

\section{Proof of Theorem~\ref{positive}}

Theorem~\ref{positive}  will be a consequence of a large deviation inequality that
holds uniformly over the whole class $\brd $ and will be proved in the
following sections.

\subsection{Time-Frequency Limiting Operators }

Let $P_R$ and $Q$ be the projection operators defined by 
\begin{equation}
  \label{eq:r2}
  P_Rf = \chi _{C_R} f \qquad \mathrm{and} \qquad Qf = \cF \inv (\chi
  _{[-1/2,1/2]^d} \hat{f} ) \, ,
\end{equation}
where $\cF^{-1}$ is the inverse Fourier transform.
Then $Q$ is the orthogonal projection from $\lrd $ onto $\cB $ and
$P_R $ is the restriction of a function to the cube $C_R$. The
composition
\begin{equation}
  \label{eq:r3}
  A_R =  Q P_R Q
\end{equation}
is the operator of time and frequency limiting. This operator has been
studied in detail by Landau, Slepian,
Pollak~\cite{LP61,LP62,SP61,slep64,sle76} and many others. It  encodes
many deep properties of \bdl\ functions and their 
restrictions. In particular, $A_R$ is a compact positive  operator of trace
class and a precisely known eigenvalue distribution. 

We summarize the properties of the spectrum that will be needed in the
sequel. 

Let $A_R ^{(1)}$ denote the operator of \tf\ limiting in dimension
$d=1$. Explicitly, $A_R ^{(1)}$  is  defined on $L^2 (\bR ) $ by the  formula
$$
(A_R^{(1)} f)\,\widehat{} \,(\xi ) = \int _{-1/2} ^{1/2} \frac{\sin \pi R (\xi
  - \eta )}{\pi (\xi -\eta)} \hat{f}(\eta) \, d\eta \qquad \text{ for
} |\xi| \leq 1/2 \, .
$$
 We denote its eigenvalues by 
$\mu _k = \mu _k  (R)$ in decreasing order and indicate the dependence on
$R$. Then  the first $[R]$ eigenvalues are 
approximately $1$,  followed by a  ``plunge region'' of thickness $\cO
(\log R)$ after which the  remaining  eigenvalues are almost zero. 
Precisely, $\mu _{[R]+1}(R) \leq 1/2 \leq \mu
_{[R]-1}(R)$; see \cite{landau93}.  This behavior of the eigenvalues is
usually formulated by saying that functions with spectrum $[-1/2,1/2]$
and ``essential'' support on $[-R/2,R/2]$  form  a
finite-dimensional subspace of ``approximate'' dimension   
$R$.  In particular, we may think of  $\brd $ as a subset of
a finite-dimensional space of dimension $R$.

The precise  asymptotic behavior of the $\mu _k$ for $k\to \infty$  was obtained by
Widom~\cite[Lemmas~1--3]{widom64}: he showed that for large $k$
\begin{equation}
  \label{eq:ch23}
  \mu _k(R) \asymp 2\pi \Big( \frac{\pi R}{8}\Big)^{2k+1} \,
  \frac{1}{k!^2} \, ,
\end{equation}
where $a_k \asymp b_k $ means that $\lim _{k\to \infty } a_k / b_k =
1$. In particular, \eqref{eq:ch23} implies  the super-exponential decay
\begin{equation}
  \label{eq:ch24}
 \mu _k(R) \leq C \exp \Big( -2k \log \big(\frac{2k}{\pi R}\big) \Big) \, .
\end{equation}
We will use the following weaker exponential estimate. 
  \begin{lemma}\cite{widom64}
    \label{widom}
Given $\alpha >0$ there exists a constant $\kappa  >0$, such that
\begin{equation}
  \label{eq:ch36}
\mu  _k (R) \leq e^{- k/\kappa } \qquad {for } \, k \geq
\frac{R}{1-\alpha} \, .
\end{equation}
  \end{lemma}
\rem\ This result is an asymptotic result for both $R\to \infty $ and $k\to
\infty $. We emphasize that the constant $\kappa $ depends only on
$\alpha $, but not on $R$.  (Widom works with the operator 
$\int _{-1} ^{1} \frac{\sin \gamma (\xi
  - \eta )}{\pi (\xi -\eta)} \hat{f}(\eta) \, d\eta $, so a simple
dilation shows that we have to use $\gamma = \pi R/2$ to obtain
$A_R^{(1)}$.)

The largest eigenvalue  $\mu _0$ of $A_R^{(1)}$  is the operator norm of $A_R^{(1)}$
and is of size $\mu _0(R) = 1 - 2 \pi \sqrt{2R} e^{-\pi R} \big( 1 +
\mathcal{O}(R\inv )\big)$ by a result of Fuchs~\cite{fuchs64}. Thus up
to terms of higher order the operator norm of $A_R$ is $\lambda _0 =
\mu _0^d = 1 - 2\pi d \sqrt{2R} e^{-\pi R}$. Assume that $\brd \neq
\emptyset$ and  $f\in \brd $, then 
$$
1-\delta \leq \int _{C_R} |f(t)|^2 \, dt = \langle A_R f ,f \rangle
\leq \lambda _0 \, .
$$
This implies that  $\delta \geq
1-\lambda _0 \geq 2\pi d \sqrt{2R} e^{-\pi R}$ (up to terms of higher order in $R$).


Next, let $C(\epsilon )$ be the   function  counting the number of eigenvalues of
$A_R ^{(1)} $ exceeding $\epsilon $, precisely
\begin{equation}
  \label{eq:ch25}
C(\epsilon ) = \mathrm{card} \{\mu _k: \mu _k  \geq
\epsilon \} \, .
\end{equation}

Then Lemma~\ref{widom}  implies that 
\begin{equation}
  \label{eq:ch26}
  C(\epsilon ) \leq \frac{R}{1-\alpha} + \kappa \log
  \frac{1}{\epsilon} \, .
\end{equation}

A different  estimate for the eigenvalue count was
obtained by Landau and Widom~\cite{LW80}: 
\begin{equation}
  \label{eq:s1}
C(\epsilon ) = R + \frac{2}{\pi} \log
  \frac{1-\epsilon}{\epsilon}  \, \log R + o \Big( \log R \Big) \, . 
\end{equation}
However, this is an asymptotic result for $R\to \infty $, and
its  proof leaves open whether  the  term $o  \big( \log R \big)$ can
be chosen independent of 
$\epsilon $. By contrast, the weaker estimate \eqref{eq:ch26} works
with a constant $\kappa $ independent of $R$, at the price of the
factor $(1-\alpha )\inv $.  Since we need the eigenvalue behavior for
\emph{fixed} $R$, we use Widom's earlier result.

Next consider the \tf\ limiting operator $A_R $ on $\lrd $
Clearly   $A_R $ is the $d$-fold tensor product of  $A_R ^{(1)}$, $A_R =
A_R^{(1)} \otimes \dots \otimes A_R^{(1)}$. Consequently,  
$\lambda $  is an eigenvalue of $A_R$, $\lambda \in \sigma (A_R)$ ,
\fif\ $\lambda = \prod 
_{j=1}^d \mu _{k_j}$, where $\mu_{k_j} \in \sigma (A_R^{(1)}) $ is an
eigenvalue of the one-dimensional 
operator $A_R ^{(1)}$. 
Since $0<\mu  _k<1 $, we have $ \prod
_{j=1}^d \mu _{k_j} \geq \epsilon$ only when $\mu _{k_j} \geq
\epsilon$ for $j=1, \dots , d$. Consequently, 
\begin{equation}
  \label{eq:s3}
  \{\lambda \in \sigma (A_R): \lambda \geq \epsilon \} \subseteq
\{\lambda = \prod _{j=1}^d \mu _{k_j}: \mu _{k_j} \in \sigma
(A_R^{(1)}),  \mu _k \geq \epsilon \}\, .
\end{equation}

We arrange the eigenvalues of $A_R$ by magnitude $1> \lambda
_1 \geq \lambda _2 \geq \lambda _3 \dots \geq \lambda _n \geq \lambda
_{n+1} \geq \dots  >0$ and again 
 let $C(\epsilon ) = \max \{n: \lambda _n \geq
\epsilon \}$ be  the function  counting the number of eigenvalues of
$A_R $ exceeding $\epsilon $. We choose $\alpha = 1/2$ and  combine
\eqref{eq:ch26} and \eqref{eq:s3}; then the eigenvalue distribution
for  $A_R$ in dimension $d$ is 
\begin{equation}
  \label{eq:s4}
  C(\epsilon ) \leq  \Big(2R + \kappa  \log
  \frac{1}{\epsilon }  \Big)^d \, ,
\end{equation}
 where $\kappa$ is independent of $R$ and $\epsilon $. 


\subsection{Covering Number  for  $\cB (R,\delta )$}

Recall that the covering numbers $N(\epsilon )= N(C,\epsilon )$
of  a compact set $C$ in a 
Banach space are defined to be the minimum number of balls of radius
less than or equal to  $\epsilon $ required to cover $C$. 
For the covering number of balls in
Euclidean space we use a   well-known estimate,
see~\cite[p.~9]{carl90} and \cite[Prop.~5]{smale2}.  

\begin{lemma}\label{euclid}
  Let $D(0,r)= \{x\in \cd : \|x\|_2 \leq r\}$ be the ball of radius
  $r$ in $\cd $. The covering number of $D(0,r)$ is given by 
  \begin{equation}
    \label{eq:s6}
    N(\epsilon ) = e^{2d  \log \frac{4r}{\epsilon} } \, .
  \end{equation}
\end{lemma}

Let us note that the covering number of the shell $D(0,r) \setminus
D(0,r(1-\delta ))$ for some $\delta >0$ is of the order $N(\epsilon )
= e^{2d  \log 
  \frac{4r}{\epsilon} } - e^{2d  \log \frac{4r(1-\delta )}{\epsilon} }  =  e^{2d  \log
  \frac{4r}{\epsilon} } \big( 1- e^{2d \log (1-\delta )}\big)$. The
difference from the covering number of  the full ball is thus negligible for large dimensions.

In the main part of the argument we will use the restriction of
bandlimited functions to the cube $C_R$. Therefore we will use the
local norms
\begin{eqnarray*}
  \|f\|_{2,R} &=& \Big( \int _{C_R} |f(x)|^2 \, dx \Big)^{1/2} \, ,\\
  \|f\|_{\infty,R} &=& \sup _{x\in C_R} |f(x)| \, ,
\end{eqnarray*}
and we denote the restriction of $\brd $ to $C_R $ by 
$$
V(R,\delta ) = P_R \brd = \{ f \in L^2(C_R): f =  \chi _{C_R} h \, \text{ for
} \, h\in \brd  \} \, .
$$

\begin{lemma}
  \label{covering}
(i) $V(R,\delta ) $ is a compact subset  in $L^2(C_R)$. 

(ii) The covering number $N_2(\epsilon )$  of $V(R,\delta) $  (with respect to $\| \cdot
\|_{2,R}$)  is bounded by 
\begin{equation}
  \label{eq:r4a}
  N_2(\epsilon ) \leq  \exp \Big( 2^{d+1}(R + \kappa  \log
  \frac{2 \sqrt{\delta}}{\epsilon }  )^d \log \frac{4\sqrt{2}}{\epsilon }\Big) \, .
\end{equation}
\end{lemma}

\begin{proof}
The finiteness of the covering numbers implies that $V(R,\delta ) $ is
compact, so it suffices to prove (ii). 

(ii) Let $\vf _n $ be the normalized  eigenfunctions of $A_R$ corresponding to the
eigenvalues of $\lambda _n$. (These are tensor products of the
standard prolate spheroidal functions.) Then $\{ \vf _n : n \in \bN \}
$ is an orthonormal basis for $\cB$. If $f = \sum _{n\in \bN } c_n \vf
_n \in \cB $, then $\|f\|_2^2 = \sum _{n\in \bN } |c_n|^2$ and
\begin{equation}
  \label{eq:rev3}
  \|f\|_{2,R}^2 = \int _{C_R} |f(x)|^2 \, dx = \langle A_Rf,f\rangle =
  \sum _{n\in \bN } |c_n|^2 \lambda _n \, .
\end{equation}
 Consequently 
$f \in \brd $ \fif\ $\mathbf{c} \in 
 S_\delta = \{ \mathbf{c} \in \ell ^2 : \|\mathbf{c}\|_2 =1 , \sum
 _{n=1}^\infty  |c_n|^2 \lambda _n \geq 1-\delta \} .$ 
Then $V(R,\delta ) $ (with the local $\| \cdot \|_{2,R}$-norm) and $S_\delta
$ (with the weighted $\ell ^2$-norm)  are isomorphic and their covering numbers
are identical.

We first determine a suitable cutoff $D$ such that  the remainder $\sum _{n>D}
|c_n|^2 \lambda _n$ is uniformly small over $S_\delta $: since
$\lambda _n <1$, we have 
\begin{eqnarray*}
  1 - \sum _{n>D} |c_n|^2 &\geq & \sum _{n\leq D} |c_n|^2 \lambda _n
  \\
&\geq & 1-\delta - \sum _{n>D} |c_n|^2 \lambda _n \\
&\geq & 1-\delta - \lambda _{D+1} \sum _{n>D} |c_n|^2 \, .
\end{eqnarray*}
We obtain that
$$
\delta \geq (1-\lambda _{D+1} ) \sum _{n>D} |c_n|^2\, ,$$
and thus
\begin{equation}
  \label{eq:rev4}
  \sum _{n>D} |c_n|^2 \lambda _n \leq \lambda _{D+1} \, \sum _{n>D}
  |c_n|^2 \leq \frac{\lambda _{D+1}}{1-\lambda _{D+1}} \delta \, .
\end{equation}
Given $\epsilon >0$, we choose the minimal $D$ so  that $\lambda _{D+1} < \frac{\epsilon
  ^2}{4\delta}$, then 
$\frac{\lambda _{D+1}}{1-\lambda _{D+1}} \delta \leq  2\lambda _{D+1}
\delta < \epsilon ^2/2 $. According to \eqref{eq:s4} we may choose $D$ to be 
\begin{equation}
  \label{eq:rev5}
  D = C\Big(\frac{\epsilon   ^2}{4\delta}\Big) \leq \Big(2R + \kappa  \log
  \frac{4\delta}{\epsilon ^2 }  \Big)^d = 2^d \Big(R+ \kappa \log
  \frac{2 \sqrt{\delta}}{\epsilon  }\Big)^d \, . 
\end{equation}

Once $D$ is determined,  choose an $\frac{\epsilon }{ \sqrt{2}}$-net
$\{\mathbf{a}_j: j=1, \dots , N_2\}$  for the unit
ball in $\bC ^D$ (with respect to the Euclidean norm) and set $f_j = \sum _{n\leq D} \mathbf{a}_j(n) \vf
_n$.  By Lemma~\ref{euclid} the cardinality of this net is 
at most $N_2=  
e^{2D \log \frac{4\sqrt{2}}{\epsilon }}$.

Given $f= \sum _{n=1}^\infty c_n \vf _n \in \brd $, choose
$\mathbf{a}_j$ and the corresponding $f_j\in \cB $, 
such that $\sum _{n\leq D} |c_n - \mathbf{a}_j(n)|^2 < \epsilon
^2/2$. Then by \eqref{eq:rev4} and the definition of $D$
\begin{eqnarray*}
  \|f-f_j\|_{2,R}^2 &=& \sum _{ n\leq D} |c_n - \mathbf{a}_j(n) |^2
  \lambda _n \, + \sum _{ n> D} |c_n|^2 \lambda _n \\
&\leq & \frac{\epsilon ^2}{2} +  \frac{\epsilon ^2}{2}  = \epsilon ^2
\, .
\end{eqnarray*}
Thus $\{ f_j \}$ is an $\epsilon $-net for $V(R,\delta )$ with respect
to $\| \cdot \|_{2,R}$. Now by Lemma~\ref{euclid} and \eqref{eq:rev5}
imply that the cardinality of this $\epsilon $-net is
 at most
 \begin{eqnarray}
N_2(\epsilon ) &\leq & \exp \Big(2D  \log \frac{4\sqrt{2}}{\epsilon} \Big) \leq  \exp \Big(2
  C\Big(\frac{\epsilon ^2}{4\delta }\Big) \log \frac{4\sqrt{2}}{\epsilon}
  \Big) \notag  \\
&\leq & \exp \Big(2^{d+1} \Big(R+ \kappa \log
  \frac{2 \sqrt{\delta}}{\epsilon  }\Big)^d \log
  \frac{4\sqrt{2}}{\epsilon}\Big) \, .    \label{eq:ch40}
  \end{eqnarray}
\end{proof}

\rem\   In dimension $d=1$  a similar estimate for the covering number
was obtained in 
~\cite{jag69}.  Estimate~\eqref{eq:ch40} also follows from general principles
in approximation theory~\cite{pinkus85} (Ch.~4.2, in particular
Thm.~2.5 and Cor.~2.6).  The estimate of the covering
number by means of the eigenvalue 
distribution, equivalently between entropy numbers and approximation
numbers, goes back  to an  
inequality of Mityagin~\cite[Ch.~9]{lorentz66}.  

As our next step we want a similar  estimate for the covering number of
$V(R,\delta )$ in the local $\| \cdot \|_{\infty ,R}$-norm.  
For this recall a   basic inequality for \bdl\ functions: \emph{If $f\in \cB $, then }
 \begin{equation}
   \label{eq:r5}
   \|f\|_\infty \leq \|f\|_2 \qquad \forall f \in \cB \, .
 \end{equation}
A similar comparison for the local norms is given in the next lemma.

\begin{lemma}
  \label{compnorm}
If $f\in \cB $, then 
\begin{equation}
  \label{eq:re1}
\|f\|_{\infty ,R}=   \max _{x\in C_R} |f(x)| \leq K_d \, 
\|f\|_2^{\frac{d}{d+2}} \, \|f\|_{2,R} ^{\frac{2}{d+2}} \, ,
\end{equation}
where the constant $K_d $ depends only on the dimension $d$ and is of
order $\cO (d)$. 
\end{lemma}

\begin{proof}
  We assume first that $f$ is real-valued and that $\alpha = \max _{x\in C_R}
  |f(x)|$ is taken at $x_0 \in C_R$. By switching sign if necessary,  we have $\alpha
  = f(x_0) \geq |f(x)|$ for all $x\in C_R$.
Next observe that by \eqref{eq:r5} and Bernstein's inequality we have 
$$
\Big\|\frac{\partial f}{\partial x_j}\Big\|_\infty \leq   \Big\|\frac{\partial
  f}{\partial x_j}\Big\|_2 \leq \pi \|f\|_2 \qquad  \text{ for } f\in \cB
\, , 
$$
and consequently
$$
\| \, |\nabla f| \, \|_\infty = \max _{x\in \rd } \Big( \sum _{j=1}^d
\Big|\frac{\partial f}{\partial x_j} (x)\Big|^2\Big)^{1/2} \leq \pi \sqrt{d}
\|f\|_2 \qquad \text{ for } f\in \cB \, .
$$

Since $f(x) = f(x_0) + \nabla f(\xi ) \cdot (x-x_0)$ for some $\xi \in
\rd $, we obtain a
lower estimate for $f$ near its maximum at $x_0$ by
$$
|f(x)| \geq \alpha - \| \, |\nabla f| \, \|_\infty \, |x-x_0| \geq
\alpha - \pi \sqrt{d} \|f\|_2 \, |x-x_0| \geq 0$$
on the ball  $B(x_0, \beta ) = \{x: |x-x_0| \leq \alpha /(\pi \sqrt{d}
\|f\|_2):= \beta \} $. We note that $\beta = \|f\|_{\infty ,R}/ (\pi
\sqrt{d} \|f\|_2) \leq (\pi \sqrt{d})\inv $ by~\eqref{eq:r5},
and thus  a fixed portion of the ball $B(x_0,\beta)$ is always
contained in $C_R$. 

Consequently (with $\sigma _{d-1}$ denoting  the surface of the
$d-1$-dimensional unit ball in $\rd $)
\begin{eqnarray*}
  \int _{C_R} |f(x)|^2 \, dx &\geq & \int _{C_R \cap B(x_0 , \beta )}
(  \alpha - \pi \sqrt{d} \|f\|_2 \, |x-x_0|)^2  \, dx
\\
&\geq & \frac{1}{2^d} \int _{B(x_0 , \beta )}
(  \alpha - \pi \sqrt{d} \|f\|_2 \, |x-x_0|)^2  \, dx
\\
&=& \frac{1}{2^d}  \int _{B(0 , \beta )} (  \alpha - \pi \sqrt{d}
\|f\|_2 \, |x|)^2  \, dx  \\ 
&=& \frac{1}{2^d} \sigma _{d-1} \, d \pi ^2 \|f\|_2^2 \,   \int _0 ^\beta
(  \beta  -  r )^2 r^{d-1}  \, dr \\
&=& \frac{1}{2^d} \sigma _{d-1} d \pi ^2\,  \|f\|_2^2  \,  
\frac{2}{d(d+1)(d+2) }\,  \beta ^{d+2} 
\end{eqnarray*}
Unraveling this inequality, we obtain that
\begin{eqnarray*}
 \|f\|_{\infty ,R} &=&  \max _{x\in C_R} |f(x)| = \pi \sqrt{d}
\|f\|_2 \,  \beta  \\
&\leq & \Big(2^{d-1} \sigma _{d-1}\inv  (d+1)(d+2)\Big)^{\frac{1}{d+2}}
  \Big( \pi  \sqrt{d} \Big) ^{\frac{d}{d+2 }} \|f\|_2 ^{\frac{d}{d+2}}    
\, \Big( \int _{C_R} |f(x)|^2 \, dx \Big) ^{\frac{1}{d+2}} \\
&=& K_d' \|f\|_2 ^{\frac{d}{d+2}}    
\, \|f\|_{2,R}^{\frac{2}{d+2}} \, . 
\end{eqnarray*}
For complex-valued $f\in \cB $, we have to take  $K_d= 2K_d' = 2
\big(2^{d-1} \sigma _{d-1}\inv  (d+1)(d+2)\big)^{\frac{1}{d+2}} 
  \big( \pi  \sqrt{d} \big) ^{\frac{d}{d+2 }}$. Using $\sigma _{d-1}
  = d\pi ^{d/2}/ \Gamma (d/2+1)$, one can then show that $K_d =  \cO (d) $.
\end{proof}

\begin{cor}
  \label{coveringinf}
(i) $V(R,\delta) $ is a compact subset in  $\cC ([-R/2, R/2]^d)$. 

(ii) The covering number $N(\epsilon )$  of $V(R,\delta) $  with
respect to $\| \cdot 
\|_{\infty,R}$  is bounded by 
\begin{equation}
  \label{eq:r4}
  N(\epsilon ) \leq  \exp \Big( 2^{d+1}\big(R + \kappa (\frac{d}{2}+1)  \log
  \frac{2 K_d}{\epsilon }  \big)^d \log
\frac{4 K_d}{\epsilon} \Big) \, . 
\end{equation}
\end{cor}

\begin{proof}
Given $\epsilon >0$, set $\epsilon _0  = 2^{-d/2} \big(\frac{\epsilon}{K_d}
\big)^{d/2+1}$ and let  $\{f_j\} $ by an $\epsilon _0$-net with respect to
$\|\cdot \|_{2,R}$. If  $f\in V(R,\delta)  $ and $\|f-f_j\|_{2,R} \leq
\epsilon _0$, then we have
$$
\|f-f_j\|_{\infty ,R} \leq K_d \|f-f_j\|_{2}^{\frac{d}{d+2}}
\, \|f-f_j\|_{2,R} ^{\frac{2}{d+2}} \leq K_d \, 2^{\frac{d}{d+2}}\,  \epsilon _0
  ^{ \frac{2}{d+2}} \leq \epsilon \, .
$$
Thus $\{ f_j\}$ is an $\epsilon $-net for $V(R,\delta )$ with respect
to $\| \cdot \|_{\infty ,R}$ and $N(\epsilon ) \leq  N_2(\epsilon _0)
$. Now use  Lemma~\ref{covering} and estimate the occurring
logarithmic term by $\log \frac{2 \sqrt{\delta}}{\epsilon_0} = \log
\frac{2 \sqrt{\delta} 2^{d/2}K_d^{d/2+1}}{\epsilon ^{d/2+1}} \leq
(\frac{d}{2}+1) \log \frac{2 K_d}{\epsilon}  $. 
  \end{proof}

The precise order of the covering number for $d=1$ with respect to the
local supremum norm $\|  \cdot \|_{\infty ,R}$ was
derived by Buslaev and Vitushkin~\cite{BV74}. Their technique is 
specifically one-dimensional and yields $N(\epsilon ) = e^{R \log
  (C/\epsilon )}$ for some constant $C>0$.

We will work with $\epsilon $-nets in the $\| \cdot \|_{\infty
  ,R}$-norm for $\epsilon = 2^{-\ell }$, $\ell
= 1 , 2, \dots $. In this case the covering number can be rewritten as
\begin{equation}
  \label{eq:s10}
    N(2^{-\ell }) \leq  \exp \Big( 2^{d+1} \Big(R+(\frac{d}{2} +1) \kappa
    \big( (\ell +1) \log 2 + \log K_d\big) \Big)^d \big(
  (\ell +2) \log 2 + \log K_d \big) \Big) 
:=  \exp p(\ell ) \, ,  
\end{equation}
where $p(\ell )=  2^{d+1} \Big(R+(\frac{d}{2} +1) \kappa
    \big( (\ell +1) \log 2 + \log K_d\big) \Big)^d \big(
  (\ell +2) \log 2 + \log K_d \big)$ is a polynomial of degree $d+1$. 

What is crucial in the above estimate, is that the  exponent grows
polynomially in $\ell $, but not faster. 

\subsection{Preparation for the proof of Theorem~\ref{positive}}

Assume that $\{x_j: j\in \bN \}$ is  an infinite sequence of i.i.d.~random
variables, each of which is uniformly distributed
over the cube  $C_R $. 

For every  $f\in \cB $ we introduce the random variable
\begin{equation}
  \label{eq:r6}
  Y_j(f)=|f(x_j)|^2- \frac{1}{R^d}\int _{C_R} |f(x)|^2 \, dx =
  |f(x_j)|^2- \E [\, |f(x_j)|^2] \, .
\end{equation}
Then $\yy $ is a sequence of independent random variables with $\E \yy
= 0$.

We first estimate  the probability distribution of the random 
  variable  
$$
\sup _{f\in \brd } \sum _{j=1}^r \yy \, .
$$
For the repeated application of Bernstein's inequality for sums of independent
random variables we will need
the following estimates for the $\yy $'s. 

\begin{lemma}\label{yyy}
Let  $f,g  \in \brd $ and $j\in \bN  $. Then 
  the following inequalities hold: 
  \begin{eqnarray}
 & & \label{eq:mo2}     \Var Y_j(f) \leq \frac{1}{R^d} \,, \\
  & & \label{eq:mo3}
  \Var (\yy - Y_j(g)) \leq \frac{4}{R^d} \|f-g\|_{\infty, R} ^2\, , \\ 
 & & \|\yy \|_{\infty } \leq 1 \, ,\\ [1ex]
 & & \label{eq:mo4}
  \|\yy - Y_j(g) \|_\infty \leq 2 \|f-g\|_{\infty , R}   \, . 
  \end{eqnarray}
\end{lemma}

\begin{proof}
We abbreviate the expected value of $|f(x_j)|^2  $ by $m (f) = R^{-d} \int
_{C_R} |f(x) |^2 \, dx $. 
Using~\eqref{eq:r5}, we obtain  
\begin{eqnarray*}
\Var Y_j (f) &=& \E [Y_j(p)^2] = \E [\,|f(x_j)|^4] - m (f) ^2 \\
&=&  \frac{1}{R^d}\int _{C_R} |f(x)|^4\, dx - m (f)^2 \\
&\leq & \frac{1}{R^d} \|f\|_{\infty , R} ^2 \, \|f\|_2^2 \leq \frac{1}{R^d} \,.
\end{eqnarray*}
 Similarly, we obtain
$$
\|Y_j(f)\|_\infty = \sup _{\omega \in \Omega } \Big|\, |f(x_j(\omega
))|^2 - m(f) \Big|  \leq \max  \Big( \|f\|_{\infty , R} ^2, \frac{1}{R^d}
\int _{C_R} |f(x)|^2 \, dx  \Big) \leq 1 \, .
$$
To prove \eqref{eq:mo3}, we write
\begin{eqnarray*}
  \Var (\yy - Y_j(g))& = &\E (\yy - Y_j(g))^2 \\
&=& \frac{1}{R^d} \int _{C_R} \big( |f(x)|^2 - |g(x)|^2 \big)^2 \, dx
- (m(f) - m(g))^2 \\
&\leq & \frac{1}{R^d} \int _{C_R}  |f(x) - g(x)| ^2 \, \big(
|f(x)|+|g(x)|\big)^2    \, dx \\
&\leq & \frac{2}{R^d} \|f-g\|_{\infty , R} ^2 \int _{\rd } (|f(x)|^2 +
|g(x)|^2 ) \, dx \leq \frac{4}{R^d} \|f-g\|_{\infty , R} ^2
\end{eqnarray*}
The last estimate follows similarly from
\begin{eqnarray*}
  \|\yy - Y_j(g) \|_\infty &\leq & \sup _{\omega \in \Omega }
\Big(  \Big| \, |f(x_j(\omega ))|^2 - |g(x_j(\omega ))|^2  \Big|  -
  \frac{1}{R^d} \int _{C_R} (|f(x)|^2 - |g(x)|^2 ) \, dx  \Big)  \\
&\leq &  \| |f|^2 - |g|^2 \|_{\infty , R}  \\
&\leq &   \|f-g\|_{\infty , R} \, \| \, |f| + |g| \|_\infty  \\
&=& 2 \|f-g \|_{\infty , R}\, .
\end{eqnarray*}  
\end{proof}

\subsection{Proof of the  sampling inequality}
\label{pr}

The sampling inequality follows from a uniform large deviation
inequality for the sampling of bandlimited functions. 

\begin{tm} \label{asym}
Let  $\{x_j : j\in \bN \}$ is a sequence of i.i.d.~random variables
that are uniformly distributed over $C_R=[-R/2,R/2]^d $.  Then there exist
constants $A,B >0$ depending on $d$ and $R$, such that
\begin{equation}
\label{eq:c11}
\bP\Bigg(\sup_{f\in \brd  } \bigg|\sum_{j=1}^r
Y_j(f)\bigg| \geq \lambda\Bigg) \leq
2A \exp\Bigg(- B\frac{\lambda^2}{41rR^{-d}+\lambda}\Bigg)
\end{equation}
for $r \in \bN $ and $\lambda \geq 0$. 

Here  $B=\frac{\sqrt{2}}{36}$. If $R$ is sufficiently large, $A$ is of
order $A=\exp (C R^d)$ for a constant depending only on $d$
and $\kappa $.   
\end{tm}

Before we prove the large deviation inequality, we show how the main
theorem follows from Theorem~\ref{asym}.

\begin{proof}[Proof of Theorem~\ref{positive}]
Choose $\lambda = \frac{r\mu}{R^d}$ and recall that $Y_j(f) =
|f(x_j)|^2 - R^{-d} \int _{C_R} |f(x)|^2 \, dx 
$. Thus the event $\cE = \{ \sup _{f\in \brd }  |\sum
_{ j=1}^r Y_j(f)| \leq r\mu R^{-d} \} $ coincides with the event 
\begin{equation}
  \label{eq:s8}
  \frac{r}{R^d} \int _{C_R} |f(x)|^2 \, dx - \frac{r\mu }{R^d} \leq
  \sum _{j=1}^r |f(x_j)|^2 \leq   \frac{r}{R^d} \int _{C_R} |f(x)|^2
  \, dx + \frac{r\mu }{R^d}  \qquad \text{ for all } f\in \brd
  \, .
\end{equation}
Since by definition $1-\delta \leq \int _{C_R} |f(x)|^2 \, dx \leq 1$, 
 we find that the event of the uniform sampling inequality
 \begin{equation}
   \label{eq:s9}
   \frac{r (1-\mu - \delta)}{R^d} \leq \sum _{j=1}^r |f(x_j)|^2 \leq
   \frac{r (1+\mu )}{R^d} \qquad \text{ for all} \,\,  f\in \brd 
 \end{equation}
is contained in $\cE $. 
As a consequence of Theorem~\ref{asym} the sampling inequality
\eqref{eq:s9} holds uniformly for all $f \in \brd $ with probability
at least $$1-2A\exp (-Br R^{-d} \mu ^2/(41+\mu )).$$


This proves Theorem \ref{positive}.
\end{proof}

We are left to prove the probability estimate of
Theorem~\ref{asym}. To estimate the probability of the deviation  of a
sum of random variables from its average  we use Bernstein's
inequality for the sums of independent random variables \cite{Ben62}: \emph{
  Let $Y_j, \jr $,  be a sequence of bounded,  independent random variables
  with $\E Y_j = 0$,  $\mathrm{Var} Y_j \leq  \sigma ^2 $, and
  $\|Y_j \|_\infty \leq M$ for $\jr $. Then} 
  \begin{equation}
    \label{eq:bern}
    \bP \Big(\Bigl|\sum _{j=1}^r Y_j\Bigr| \geq \lambda \Big) \leq 2
    \exp \Big( -     \frac{\lambda ^2}{2 r \sigma ^2 + \frac{2}{3}
      M\lambda} \Big) \, .  
  \end{equation}

\begin{proof}[Proof of Theorem~\ref{asym}]
\textit{Step {\rm 1:} A metric entropy argument.}
 For a given 
$\ell  \in \bN$, we construct an $2^{-\ell} $-covering  for $
V(R,\delta ) $
with respect to the local norm
$\| \cdot \|_{\infty ,R}$. Let $\cA _\ell $ be the corresponding $2^{-\ell
}$-net for $\ell = 1, 2, \dots $. 
Then $\cA _{\ell} $ has cardinality at most $N(2^{-\ell}) \leq
e^{p(\ell )}$ for some polynomial of degree $d+1$ by
Corollary~\ref{coveringinf}. 


Given  $f\in \brd$, let $f_\ell $ be the function in $ \cA (2^{-\ell
})$ that is closest 
to $f$ in $\|\cdot \|_{\infty ,R}$-norm, with some convention for breaking
ties.
Since $\|f-f_{\ell }\|_{\infty ,R} \to 0$ 
 we  can write 
$$
Y_j(f)= Y_j(f_1)  +(Y_j(f_2)-Y_j(f_1))+(Y_j(f_3)-Y_j(f_2))+ \cdots.
$$
If $\sup_{f\in \brd } |\sum_{j=1}^r Y_j(f)| \geq 
\lambda$, then $ \cE_\ell$ must hold for some $\ell\geq 1$, where
 $$ \cE_1=\Big\{\mbox{there exists }f_1 \in \cA (1/2)\mbox{ such that }|\sum_{j=1}^r Y_j(f_1)|
\geq \lambda /2\Big\}$$ and 
\begin{eqnarray*}
\cE_\ell&= &\Big\{ \text{there exist}  \,   f_\ell  \in \cA(2^{-\ell}), \, \text{
  and } f_{\ell -1} \in \cA(2^{-\ell +1})\, \text{ with  }\notag\\
& &\qquad \norm{f_\ell-f_{\ell -1}}_{\infty, R}\leq
3 \cdot 2^{-\ell}, \notag  \\
 &  & \qquad\text{ such that }  \Big|\sum_{j=1}^r \big(
Y_j(f_\ell)-Y_j(f_{\ell -1})\big)\Big|\geq\lambda/2\ell ^2\Big\} \, .   \label{eq:ch37}
\end{eqnarray*}
If this were not the case, then, with $f_0=0$,  
\begin{equation*}
  \Bigg|\sum_{j=1}^r Y_j(f)\Bigg| \leq  \sum
_{\ell=1}^\infty \Bigg| \sum_{j=1}^r (Y_j(f_\ell) -
Y_j(f_{\ell-1})) \Bigg|  
\leq  \sum _{\ell=1}^\infty \frac{\lambda }{2\ell^2}
= \frac{\pi ^2}{12} \lambda < \lambda. 
\end{equation*}

 Next  we
estimate the probability of $\cE _\ell $.

\textit{Step $2$.} 
We estimate the term $\ell =1 $ separately. 
 For fixed $f\in \cA(1/2)$, the probability of the event $\cE _1$   is
bounded, using Bernstein's inequality~\eqref{eq:bern} and
Lemma~\ref{yyy}, by 
\begin{equation*}
  2\exp\bigg(-\frac{\lambda^2/4}{2r\Var Y_j(f)+\frac23 (\lambda/2)
\norm{Y_j(f)}_\infty}\bigg) \leq 2 \exp \bigg( - \frac{\lambda ^2}{ 2r R^{-d} +
 \lambda /3} \bigg) .
\end{equation*}
There are at most $N(1/2) = \exp \Big( 2^{d+1}(R + \kappa (\frac{d}{2}+1)  \log
  4 K_d  )^d \log 8 K_d  \Big)  $
functions  in $\cA(1/2)$, so the probability of $\cE    _1$  is bounded by 
\begin{equation}
\label{eq:c6}
2  \exp \Big( 2^{d+1}(R + \kappa (\frac{d}{2}+1)  \log
  4 K_d  )^d \log 8 K_d  \Big)  \,  \exp \bigg( - \frac{\lambda ^2}{ 2r R^{-d} +
 \lambda /3} \bigg) \, .
\end{equation}

\textit{Step $3$.} For $\ell \geq 2$, we estimate the probability of
$\cE _\ell $  in a similar fashion by using
Lemma~\ref{yyy}, \eqref{eq:mo3}, and~\eqref{eq:mo4}. If 
$f\in \cA(2^{-\ell}) $ and $g\in \cA(2^{-\ell +1})$ with
$\norm{f-g}_{\infty, R}\leq 3\cdot 2^{-\ell}$, 
we have
\begin{eqnarray*}
\bP\bigg(  \bigg|\sum_{j=1}^r (Y_j(f)&-&Y_j(g))\bigg|
>\frac{\lambda}{2\ell ^2}\bigg)\\ 
&\leq & 2\exp\bigg(-\frac{\lambda^2/4\ell ^4}{2r \cdot 4\cdot R^{-d}
   (3\cdot 2^{-\ell})^2 +\frac 23 2 \cdot 3 \cdot 2^{-\ell-1} 
\lambda/\ell ^2}\bigg) \\
&=& 2 \exp \Big( - \frac{2^\ell}{8 \ell ^2} \, \frac{\lambda
  ^2}{36rR^{-d} \ell ^2 2^{-\ell } + \lambda} \Big) \, .
\end{eqnarray*}

Note $36 \ell^2/2^\ell< 41 $.
There are at most $N(2^{-\ell })$ functions in $\cA (2^{-\ell})$
and $N(2^{-\ell +1})$ functions in $\cA (2^{-\ell+1})$. 
Finally, this can happen for
any $\ell$. So 
the probability of $\bigcup _{\ell =2}^\infty \cE _\ell $  is bounded by
\begin{eqnarray}
  \label{eq:ds1}
 \sum_{\ell=2}^\infty & & 
  N(2^{-\ell }) N(2^{-\ell +1}) 2\exp \Big( - \frac{2^\ell}{8 \ell
    ^2} \, \frac{\lambda 
   ^2}{41rR^{-d}  + \lambda} \Big)  \\
&\leq & \sum_{\ell=2}^\infty  2\exp \Big( p(\ell ) + p(\ell -1)  -
\frac{2^\ell}{8 \ell     ^2} \, \frac{\lambda 
   ^2}{41rR^{-d}  + \lambda} \Big) \, , \notag 
\end{eqnarray}
where we use~\eqref{eq:s10} for the covering number.

\textit{Step $4$.}
We will need the following inequality:            

\medskip
\textit{If $p , a >0$, then 
\begin{equation}\label{E32a}
\sum _{\ell = 2} ^\infty e^{-a^\ell p } \leq \frac{1}{p a  \log a}
e^{-ap } \, .
\end{equation}}

  This inequality follows from the integral test and the  substitution $a^x=u$:
  \begin{eqnarray*}
    \sum _{\ell = 2} ^\infty e^{-a^\ell p } &\leq & \int
    _{1}^\infty e^{-a^xp } dx \\
&=& \frac{1}{\log a}  \int _{a} ^\infty e^{-p u} \, \frac{du}{u}\\
&\leq & \frac{1}{a \log a}  \int _{a} ^\infty e^{-p u} \, du \\
&=& \frac{1}{p a  \log a} e^{-ap } \, .
  \end{eqnarray*}

\textit{Step $5$.} 
To estimate the sum~\eqref{eq:ds1}, we rewrite and simplify each term.
Set 
\begin{eqnarray}
\psi &=& \frac{\lambda^2}{41 r R^{-d} + 
  \lambda } \label{eq:ch38a}\\  
    c_1 & = &\min _{\ell \geq 2} \frac{2^{\ell /2}}{8\ell ^2} \label{eq:ch39} \\
c_2 &=& \max _{\ell \geq 2} \frac{2p(\ell )}{2^{\ell /2}} \, .   \label{eq:ch38}
\end{eqnarray}
and $\sup _{\ell \geq 2} \ell ^2 /
2^\ell = 9/8$. 
Then the $\ell$-th term in \eqref{eq:ds1} is majorized by 
$$
\exp \Big( -2^{\ell /2} ( c_1 \psi  -c_2) \Big)  \, .
$$
If $\psi >0$ is large enough so that $p:= c_1 \psi  -c_2 >0$, then (\ref{E32a})
implies that 
\begin{eqnarray}
\P (\bigcup _{\ell = 2 } ^\infty \cE _\ell ) &\leq & 2 \frac{1}{(c_1\psi - c_2)
  \sqrt{2}  \log \sqrt{2}} \, e^{-\sqrt{2} (c_1\psi -c_2) }  \notag \\ & = &
\frac{2\sqrt{2}}{\log 2}\,\frac{e^{\sqrt{2}c_2}}{c_1 \psi - c_2} \,
  \exp \Big(- \frac{\sqrt{2} c_1\lambda^2}{41 r R^{-d} +\lambda } \Big)
  \,  . \label{eq:c6a}
\end{eqnarray}

Since the term for $\ell = 1$ has the same form, we have proved that 
$$\bP\Bigg(\sup_{f\in \brd  } \bigg|\sum_{j=1}^r
Y_j(f)\bigg| \geq \lambda\Bigg) \leq
2A \exp\Bigg(- \frac{\sqrt{2}c_1\lambda^2}{41r R^{-d}+\lambda}\Bigg) \, ,
$$
whenever $\psi > c_2/c_1$.

For the exponent $B$ we may take the smaller of the exponents in
\eqref{eq:c6} and \eqref{eq:c6a}, i.e., $B=\min (3, \sqrt{2}c_1)$. If we choose
$\lambda $ large enough, so that $c_1 \psi -c_2 \geq
\frac{2 \sqrt{2}}{\log 2}$, then we may
take $A= \max ( \exp \Big( 2^{d+1}(R + \kappa (\frac{d}{2}+1)  \log
  4 K_d  )^d \log 8 K_d  \Big),
e^{\sqrt{2} c_2})$. Thus we have proved Theorem~\ref{asym}.

\emph{Step 6.} To obtain an idea of the magnitude of the constants
involved, we give some rough estimates for $c_1$ and $c_2$, $A$
and $B$.  

For $c_1$ we obtain 
$$
c_1 = \frac 18\,  \min _{\ell \geq 2} \frac{2^{\ell /2}}{\ell ^2} =
\frac{1}{36} \, , 
$$
so the exponent $B$ in~\eqref{eq:c11} is  $\sqrt{2} c_1
= \frac{\sqrt{2}}{36}$, which is approximately $ \approx 0.0393 $.


As for $c_2$, recall that $p(\ell )=  2^{d+1} \Big(R+(\frac{d}{2} +1) \kappa
    \big( (\ell +1) \log 2 + \log K_d\big) \Big)^d \big(
  (\ell +2) \log 2 + \log K_d \big)$. If $(\frac{d}{2} +1) \kappa
    \big( (\ell +1) \log 2 + \log K_d \leq R$, then  
$$
\frac{p(\ell )}{ 2^{\ell / 2}} \leq 2^{d+1} (2R)^d  \max _{\ell \geq
  2} \frac{  (\ell +2) \log 2 + \log K_d}{2^{\ell / 2}} \leq c_3
R^d \, .
$$ 
In the other case, we may estimate $p(\ell ) / 2^{\ell /2}$ by a
constant that depends on $d$, $K_d$ and $\kappa $, but not on
$R$. Thus for $R$ sufficiently large, we obtain $c_2 \leq c_3
R^d$ and $A \leq \exp ( C R^d)$.   

Finally consider the condition  $c_1\psi - c_2 \geq
\frac{2^{3/2}}{\log 2}$, which follows from 
$\frac{ \lambda^2}{41 r R^{-d} +\lambda 
} \geq c_4 R^d \geq  \frac{c_2+2^{3/2}/\log 2}{c_1}$. Since
$x\geq B+\sqrt{D} $ implies 
$x^2 \geq Bx+D$, we find that 
$$\lambda \geq c_4 +(41c_4r R^{-d})^{1/2}\, ,
$$ 
for a constant independent of $R$. 
\end{proof}

\def\cprime{$'$} \def\cprime{$'$}


\end{document}